\documentclass[11pt]{amsart}
\usepackage{amssymb}

\newtheorem{thm}{Theorem}

\newtheorem{prop}[thm]{Proposition}
\newtheorem{lemma}[thm]{Lemma}

\theoremstyle{remark}

\theoremstyle{definition}
\newtheorem{defn}[thm]{Definition}

\numberwithin{equation}{section}

\newcommand{\R}{{\mathbb R}}
\newcommand{\C}{{\mathbb C}}

\newcommand{\Q}{{\mathbb Q}}

\renewcommand{\P}{{\mathbb P}}
\newcommand{\sgn}{\operatorname{sgn}}

\newcommand{\Sym}{\operatorname{Sym}}
\newcommand{\Span}{\operatorname{Span}}
\newcommand{\Res}{\operatorname{Res}}

\newenvironment{pf}{\smallskip\noindent{\bf Proof.}\ }{\qed\smallskip}

\title{Weak approximation for linear systems of quadrics}

\author{Bo-Hae Im}
\address{Department of Mathematics, University of Utah,
Salt Lake City, UT 84112-0090, USA}
\email{im@math.utah.edu}

\author{Michael Larsen}
\address{Department of Mathematics, Indiana University,
Bloomington, IN 47405, USA}
\email{larsen@math.indiana.edu}

\thanks{The second named author was partially supported by NSF grant DMS-0100537. }

\begin{document}

\begin{abstract}
We give local conditions at $\infty$ ensuring that the intersection of $n$ quadrics in $\P^N$, $N\gg n$, satisfies weak approximation.
\end{abstract}

\maketitle

Let $K$ be a number field, $V$ a finite-dimensional vector space over $K$, and
$W\subset \Sym^2 V^*$ a $K$-vector space of quadratic forms on $V$.   We can regard $W$ as a linear system of quadrics in $\P V$, and we denote the intersection of all quadrics in the system $X_W$.    
If $S$ is a finite set of places of $K$, we can ask whether $X_W$ satisfies weak approximation with respect to $S$, i.e., whether the diagonal map
$$X_W(K)\to\prod_{s\in S} X_W(K_s)$$
has dense image.   Assuming that $X_W(K_s)\neq\emptyset$ for all $s\in S$, this is a refinement of the basic question of whether $X_W(K)$ is non-empty.

For $\dim W = 1$, the Hasse-Minkowski theorem implies $X_W(K)$ is non-empty if $\dim V\ge 5$ and a generator of $W$ is indefinite at every real place of $K$.
The weak approximation theorem for quadric hypersurfaces \cite{PR}~\S7.1~Cor.~1 then applies.
For $\dim W = 2$, 
J-L.~Colliot-Th\'el\`ene, J-J.~Sansuc, and P.~Swinnerton-Dyer  \cite[10.4]{CSS2} gave an affirmative answer to a question of A.~Pfister \cite{Pf}: if $\dim V > 4\dim W$ and $K$ is totally imaginary (Pfister originally raised the question in the case $K=\Q(i)$), must $X_W(K)$ be non-empty?   D.~Leep \cite[2.7]{Le} proved a qualitative version of Pfister's conjecture: if 
$$\dim V \ge 2(\dim W)^2+2\dim W + 1,$$
and $K$ is totally imaginary, then $X_W(K)$ is non-empty.  Using the circle method,
W.~Schmidt \cite{Sc} proved an analogous result over $\Q$, under the hypothesis that 
$X_W(\R)$ has a non-singular point.

Weak approximation is more delicate.  A simple relation between the number of variables and
the number of equations is no longer possible, even in the case of totally imaginary fields.
For instance, if $\alpha\in K$ does not have a square root in $K$, then the quadric hypersurface
given by $x_1^2-\alpha x_2^2=0$ fails to satisfy weak approximation no matter how large its dimension.
The obstruction can be understood in terms of the rank of the defining quadratic form
being too low.
For $W$ generated by a single indefinite quadratic form of rank $\ge 5$, weak approximation of $X_W$ is well known (see, e.g., \cite{PR}~S7.1~Cor.~1.)  Colliot-Th\'el\`ene, Sansuc, and Swinnerton-Dyer have a weak approximation theorem (\cite{CSS1}, \cite{CSS2}) for $\dim W = 2$, $\dim V\ge 9$,
assuming every form is indefinite and has rank $\ge 5$.
See also \cite{SS} and \cite{Sko} for good bounds when $\dim W\le 3$.

The main result of this paper is a weak approximation theorem for $X_W$ for arbitrary $\dim W$
and arbitrary number fields $K$.
The hypotheses make no explicit demands concerning the singular locus of $X_W$; instead, we impose conditions on the rank and signature of non-zero forms in $W$ with respect to the archimedean places of $K$.

There seems to be a general philosophy that a complete intersection $X$ of $n$ hypersurfaces
of degree $\le d$ should have many rational points as long as there are no local obstructions (especially no obstruction at real places) and $\dim X$ is large compared to $n$, $d$, and the dimension of the singular locus of $X$.
This paper exhibits one concrete instance of this philosophy for $d=2$.
For general $d$, C.~Skinner \cite{Ski} has proven weak approximation for smooth complete intersections as long as a certain auxiliary variety is not too badly singular.  In particular, he has proven weak approximation for non-singular cubic hypersurfaces of sufficiently high dimension.
It is not unlikely that his work can provide an alternative route to our result.

We would like to express our appreciation of the help we received from 
Colliot-Th\'el\`ene, whose pointed questions
and suggestions led to considerable improvements in the exposition of this paper.  
The original admissibility bound in our theorem was not stated 
(and was in fact worse than exponential in $\dim W$).
Both Colliot-Th\'el\`ene and the referee strongly urged us to make it explicit, and thereby encouraged
us to improve it to its current form (which is quadratic).  We are grateful to both of them.

\bigskip

We begin with a definition intended at the same time to control the singular locus and to guarantee that there is no real obstruction:

\begin{defn}
If $V$ is a vector space over $k=\R$ or $k=\C$, $W\subset \Sym^2 V^*$ a vector space of quadratic forms on $V$, and $m$ is a positive integer, we say 
$W$ is \emph{$m$-admissible} if for every non-zero $w\in W$, there are at least $m$ hyperbolics contained in $w$.  When $k=\C$, this just means that the rank of $w$ is at least $2m$; over $\R$ it means $V$ contains $m$-dimensional subspaces $V_w^+$ and $V_w^-$ on which 
$w$ is positive-definite and negative-definite respectively.  If $V$ is a vector space over a
number field $K$, $W\subset \Sym^2 V^*$, and $m$ is a positive integer, we say
$W$ is \emph{$m$-admissible} if $W\otimes 1$ is an $m$-admissible space of quadratic
forms on $V\otimes_K K_s$ for every archimedean completion $k=K_s$ of $K$.
\end{defn}

We can now state the main theorem.

\begin{thm}
\label{main}
Let $f(n) = 2n^2+2n-1$.
For every positive integer $n$, number field $K$, $K$-vector space $V$, and
$n$-dimensional vector space of quadratic forms $W\subset \Sym^2 V^*$
which is $f(n)$-admissible, the variety $X_W$ always satisfies weak approximation over $K$.
That is, for every finite set $S$ of places of $K$,
$X_W(K)$ is dense in $\prod_{v\in S} X_W(K_s)$.
\end{thm}

The growth rate of $f(n)$ is quadratic like those of \cite{Le} and \cite{Sc}, but the actual minimum
$\dim V$ to which the theorem can be applied is larger than the bounds of
Leep and Schmidt.
We have $f(2) = 11$, $f(3) = 23$, so our result certainly does not improve existing bounds in
low dimensions.  The virtues of our theorem are that it gives full weak approximation rather than mere existence of solutions, and it does so for an arbitrarily large system of
quadratic forms and without explicit hypotheses on the geometry of $X_W$.  
Moreover, the proof is elementary.
To find points in $X_W(K)$,  we choose a quadratic form $w_0\in W$ and a complementary subspace 
$W^0\subset W$.  We then proceed to find five vectors $v_1,\ldots,v_5$ which are
mutually orthogonal with respect to $W$ and isotropic with respect to $W^0$
but not with respect to $w_0$.    Then we apply weak approximation to the locus $w_0=0$
on the vector space spanned by the $v_i$.

The following notation will be useful.
Given vector spaces $V$,  $W\subset\Sym^2 V^*$, and $V_0\subset V$, we have a natural surjective linear transformation $\Sym^2 V^*\to \Sym^2 V_0^*$.  The image of $W$ under this transformation will be denoted $\Res_{V_0} W$.  Note that if $\ker V^*\to V_0^*$ is spanned by linear forms
$L_1,\ldots,L_r$ on $V$, a quadratic form on $V$ lies in $\ker\Sym^2 V^*\to \Sym^2 V_0^*$ 
if and only if it is of the form $\sum_{i=1}^r L_i M_i$ for $M_i\in V^*$.  In particular, the kernel of
restriction consists only of forms of rank $\le 2r$.

\begin{lemma}
\label{restrict}
Let $V$ be a finite-dimensional vector space over $k=\R$ or $k=\C$ 
and $W\subset \Sym^2 V^*$ be a vector space of quadratic forms.  Let $V_0\subset V$ be a subspace of codimension $d$.   If $W$ is $m$-admissible for some integer $m>d$, then $\Res_{V_0}W$ is $m-d$-admissible.
\end{lemma}

\begin{proof}
Let $\bar w\in \Res_{V_0}W$ be non-zero, and let $w$ be an element of $W$ mapping to
$\bar w$.  If $k=\C$, there there exists a $2m$-dimensional subspace $V_w\subset V$
in which $w$ is non-degenerate.  Then $V_w\cap V_0\subset V_0$ is a 
subspace of dimension $\ge 2m-d > 2(m-d)$ on which $\bar w$ is non-degenerate.
For $k=\R$, let $V_w^+$ (resp. $V_w^-$) denote an $m$-dimensional
subspace of $V$ on which $w$ is positive-definite (resp. negative-definite).  Then
$\bar w$ is positive-definite (resp. negative-definite) on $V_w\cap V_0$, which has dimension
$\ge m-d$.
\end{proof}

Given a real or complex vector space $V$ and a vector space $W\subset \Sym^2 V^*$
of quadratic forms, we define the (non-linear) \emph{evaluation map} $E:\; V\to W^*$ by the
identity
$$E(v)(w) = w(v).$$
\begin{prop}
\label{sphere}
Let $g(n)=n^2-n+1$.
Let $n$ be a positive integer, $V$ a finite-dimensional vector space over $k=\R$ or $k=\C$, and 
$W\subset \Sym^2 V^*$ an 
$n$-dimensional space of quadratic forms.  If $W$ is $g(n)$-admissible, then the evaluation map $E:\;V\to W^*$ is surjective.
\end{prop}

\begin{proof}
The idea of the proof is that for fixed $v_1,v_2\in V$, the map 
$$(a,b)\mapsto E(av_1,bv_2)$$
is of the form 
$$(a,b)\mapsto a^2 w_1^* + b^2 w_2^* + ab w_{1,2}^*$$
for suitable 
vectors $w_1^*=E(v_1)$, $w_2^*=E(v_2)$, and $w_{1,2}^*\in W^*$.  
In general, the image is therefore a cone in
a 3-dimensional subspace of $W^*$.  If, however, it happens that $w_{1,2}^*$ is a linear
combination of $w_1^*$ and $w_2^*$, then this formula can be used to show that
a wide range of linear combinations of $w_1^*$ and $w_2^*$ belong to $E(V^*)$.
For general $v_1$ and $v_2$, there is no reason to believe this condition will hold, but
we can use an inductive argument to prove that there exists $v_3\in V$ such that
the analogous conditions hold for $v_1$ and $v_3$ and for $v_2$ and $v_3$, and moreover
the vector $w_3^* = E(v_3^*)$ belongs to the space of $w_1^*$ and $w_2^*$.  The details are somewhat simpler in the case $k=\C$, so we begin there.

Suppose $k=\C$.
It suffices to prove that $E(V)$ is a vector space, since if $w(E(V))$ = 0 for some vector $w\in W$,
then $w(v)$ = 0 for all $v\in V$, and this implies $w=0$.  Let $v_1,v_2\in V$, and let
$w_i^* = E(v_i)$.  We would like to show that every linear combination of $w_1^*$ and $w_2^*$
lies in $E(V)$.  As $E(av) = a^2 E(v)$, this is obvious unless $w_1^*$ and $w_2^*$ span
a two-dimensional space.  

Let $W_1 = (W^*/\C w_1^*)^*$ denote the hyperplane in $W$ determined by $w_1^*$.
Let 
$$V_{1,2} = \{v\in V\mid w(v+v_1) = w(v)+w(v_1),\;w(v+v_2)=w(v)+w(v_2)\ \forall w\in W_1\}$$ 
be the subspace of $V$ consisting of vectors orthogonal to $v_1$ and $v_2$ with respect to every
quadratic form in $W_1$ (i.e., with respect to every symmetric bilinear form associated to a quadratic form in $W_1$).  It is enough to check this orthogonality on a basis for $W_1$, so 
$V_{1,2}$ is of  codimension $\le 2n-2$ in $V$.  By Lemma~\ref{restrict}, $\Res_{V_{1,2}}W_1$
is $g(n)-(2n-2) = g(n-1)$-admissible.  
Note that as $g(n) \ge 2n-2$, the restriction map of $W_1$ to $V_{1,2}$ is an isomorphism.
By the induction hypothesis, there exists $v_3\in V_{1,2}$
whose image in 
$$\Res_{V_{1,2} }W_1^* = W_1^* = W^*/\C w_1^*$$
equals the class of $w_2^*$ in $W^*/\C w_1^*$.  Let $w_3^* = E(v_3)$.  
Thus $\Span(w_1^*, w_3^*) = \Span(w_1^*,w_2^*)$
and $v_1$ and $v_3$ are orthogonal with respect to every form in $W_1$.  This last condition implies
that there exists $c\in \C$ such that
$$w(av_1+bv_3) = a^2 w(v_1) + b^2 w(v_3) + abc w_1^*(w) = ((a^2+abc)w_1^* + b^2 w_2^*)(w).$$
For any $c\in\C$, the function $(a,b)\mapsto (a^2+abc,b^2)$ is a surjective self-map of $\C^2$,
so indeed $E$ is surjective.

Now suppose $k=\R$.  The formula $E(av) = a^2 E(v)$ implies $E(V)$ is a cone.  If $E(V)$ is contained in a closed half-space $H_w = \{w^*\mid w(w^*)\ge 0\}$
then $w$ is semidefinite, contrary to hypothesis.  If $E(V)$ is convex, then
by \cite[11.7.3]{Ro}, it follows that $E(V) = W^*$.  To prove convexity, it suffices to prove
that the line segment connecting any $w_1^* = E(v_1)$ and $w_2^* = E(v_2)$ lies in $E(V)$.
This is clear if $w_1^*$ and $w_2^*$ are linearly dependent, so we may assume
that $\Span(w_1^*,w_2^*)$ is a plane.  As $E(V)$ is a cone, it suffices to prove that there is a
path connecting $w_1^*$ and $w_2^*$ which lies entirely in a closed half-plane within
that plane.

We construct $V_{1,2}$, $v_3$, and $w_3^*$ as in the complex case.  We define $c,d\in\R$ such that
\begin{gather*}
w(av_1+bv_3) = a^2 w(v_1) + b^2 w(v_3) + abc w_1^*(w)\\ 
w(av_2+bv_3) = a^2 w(v_2) + b^2 w(v_3) + abd w_1^*(w)\\
\end{gather*}
Replacing $v_1$ and $v_2$ by $-v_1$ and $-v_2$ respectively if necessary, we may assume
$c,d\ge 0$.  By construction, $w_1^*$, $w_2^*$, and $w_3^*$ lie in a common half-plane.
It follows that 
\begin{gather*}
E(av_1+bv_3) = a^2 w_1^* + b^2 w_3^* + abc w_1^*\\ 
E(av_2+bv_3) = a^2 w_2^* + b^2 w_3^* + abd w_1^*\\
\end{gather*}
all lie in this half-plane.  Specializing to $(a,b) = (t,1-t)$ for $t\in[0,1]$, we conclude that
there exist paths from $v_3$ to $v_1$ and from $v_3$ to $v_2$ contained in a common
half plane and therefore that there exists a path from $v_1$ to $v_2$ in that half-plane.  The proposition follows.

\end{proof}

We can now prove the main theorem.

\begin{pf}
We use induction on $n$.  The case $n = 1$ is well-known, so we may assume $n\ge 2$.  We fix a non-zero $w_0\in W$ and a codimension-1 subspace $W^0\subset W$ such that $w_0\not\in W^0$.
Let $S$ be a finite set of places of $K$ and $S_\infty$ be the set of archimedean places of $K$.
For each $s\in S$, we fix a non-empty open set $U_s$
in $\P (V\otimes_K K_s)$.
We apply the induction hypothesis to $W^0$ and $\prod_{s\in S} U_s\cap X_{W^0}(K_s)$
to obtain a vector $v_1\in V$ such that
$v_1$ is isotropic for all quadratic forms in $W^0$ and 
for all $s\in S$, the line $K_s v_1$ defines a point in $\P (V\otimes_K K_s)$ which lies in $U_s$.
If $w_0(v_1)=0$, we are done, so we assume this quantity is not zero.
Next we construct vectors $v_2,\ldots,v_5\in V$ iteratively with the following properties:
\begin{enumerate}
\item For $i > j \ge 1$, the vectors $v_i$ and $v_j$ are mutually orthogonal with respect to every form
$w\in W$.
\item For $i\ge 2$, $v_i$ is $W^0$-isotropic
\item For $i\ge 2$, $v_i$ is not $w_0$-isotropic.
\item For every real embedding $s$ and every $i\ge  2$, we have
$\sgn s(w_0(v_i)) = - \sgn s(w_0(v_1)).$
\end{enumerate}

To achieve these properties, define for $2\le i\le 5$
$$V_i = \{v\in V\mid w(v+v_j) = w(v)+w(v_j)\ \forall j < i\}.$$
Clearly $V_i$ is of codimension
$\le (i-1)n\le 4n$ in $V$.  By Lemma~\ref{restrict}, $\Res_{V_i}W$ is
$f(n)-4n$-admissible, and as $f(n)-4n = f(n-1)\ge n^2-n+1$, Proposition~\ref{sphere} applies.
As the rank of every non-zero $w\in W$ is $\ge 2g(n) > 4n$, the restriction map
$W\to \Res_{V_i}W$ is an isomorphism.
This means that for every archimedean completion $K_s$ of $K$, there exists 
$v_{i,s}\in V_i\otimes_K K_s$ which is isotropic for every quadratic form in $W^0$ and satisfies
$w_0(v_{i,s})\neq 0$.  If $K_s=\R$, we can additionally choose the sign of $w_0(v_{i,s})$ to be 
$-\sgn s(w_0(v_1))$.   For each archimedean place $s$ of $K$, we define a non-empty neighborhood 
$U_{i,s}$ of $K_s v_{i,s}$
in $X_{\Res_{V_i}W^0}(K_s)$  given by
\begin{equation*}
U_{i,s} = 
\begin{cases}
\{K_s v\mid v\in V_i, w_0(v)\neq 0\}& \textrm{if $K_s=\C$} \\
\{K_s v\mid v\in V_i, s(w_0(v) w_0(v_1)) < 0\}& \textrm{if $K_s = \R$}.\\
\end{cases}
\end{equation*}
We now apply the induction hypothesis to $\Res_{V_i}W^0$ and $\prod_{s\in S_\infty} U_{i,s}$
to find the vector
$v_i$.  Condition (1) is automatic since $v_i\in V_i$.   Condition (2) holds
because $K v_i$ lies in $X_{\Res_{V_i}W^0}(K)$.  Conditions (3) and (4) are true because
$K v_i$ lies in $\prod_{s\in S_\infty} U_{i,s}$.

Now consider the restriction of $w_0$ to $V_0 := \Span(v_1,v_2,v_3,v_4,v_5)$.
We know that $w_0$ is diagonal in terms of the basis $v_i$, of full rank, and
indefinite at all real places of $K$ by conditions (1), (3), and (4) respectively.  
Moreover, the intersection
$$P V_0 \cap \prod_{s\in S} U_s$$
is open and non-empty (because it contains the line $K v_1$).  By the weak approximation
theorem for indefinite quadratic forms of rank 5, there exists $v\in V_0$ which is
isotropic for $w_0$ and such that $K v \in \prod_{s\in S} U_s$.
Conditions (1) and (2) above now imply that $v$ is also isotropic
for $W^0$ and therefore for all of $W$.
The theorem follows.

\end{pf}

We conclude with several questions.   Does the theorem remain true if we replace our notion of $m$-admissibility by the condition that every non-zero $w\in W$ is indefinite and non-degenerate on an $m$-dimensional subspace of $V$?  
Do we get weak approximation if we assume only that $\dim V\gg 0$ and every non-zero $w\in W$
is indefinite and of rank $\ge 5$?

Our proof of Proposition~\ref{sphere} depends on a trick that shows that under suitable admissibility conditions, the real locus of a system of quadratic equations (in real projective space) is connected.   Can this trick be used to analyze higher connectivity?  For instance, if $V$ is a real vector space, $k$ is a positive integer, and $W$ a vector space of quadratic forms which is $m$-admissible for some $m$ depending on $k$ and $\dim W$, can we conclude that $X_W(\R)$ has a $k$-connected double cover?

\end{document}